\documentclass[11pt,a4paper]{article}
\usepackage{a4wide}
\usepackage[utf8]{inputenc}
\usepackage{amsmath}
\usepackage{amsfonts}
\usepackage{amssymb}
\usepackage{graphicx}
\usepackage{caption}
\usepackage{fancyhdr}
\usepackage{fleqn}
\usepackage{latexsym}
\usepackage{float}
\usepackage[bookmarks=false]{hyperref}
\usepackage{xcolor}
\usepackage{enumitem}
\setlist{nosep}
\author{Benne de Weger \\ Department of Mathematics and Computer Science \\ 
       Eindhoven University of Technology \\ \texttt{b.m.m.d.weger@tue.nl}}
\newcommand{\titel}{A lot of fudge around $ A + B = C $}
\newcommand{\kortetitel}{\titel}
\newcommand{\version}{version 1.0, \today}
\title{\titel}
\date{\version}

\newtheorem{lemma}{Lemma}
\newtheorem{conjecture}[lemma]{Conjecture}

\newcommand{\Sha}{\text{Sha}}
\newcommand{\Q}{\mathbb{Q}}

\begin{document}

\maketitle

\begin{abstract}
\noindent This paper reports on experiments searching for elliptic curves over $ \Q $ with large 
Tamagawa products. The main idea is to look among curves related to good $ abc $-triples.
\end{abstract}

\section{Introduction}
Let $ E $ be an elliptic curve defined over $ \Q $, with conductor $ N = N(E) $ and Tamagawa product 
(also called `fudge factor') $ \tau = \tau(E) $. In my 1998 paper \cite{dW2} I proved a (conditional) 
upper bound for $ \tau $ in terms of $ N $, namely 
\begin{lemma}[de Weger, 1998]
Szpiro's Conjecture implies that for all elliptic curves 
\[ \tau \ll N^{\epsilon} . \]
\end{lemma} 
This is short for: For each $ \epsilon > 0 $ there exists an absolute constant $ C_{\epsilon} $ such
that for all elliptic curves over $ \Q $ it is true that $ \tau < C_{\epsilon} N^{\epsilon} $.
In fact, I proved a slightly better result, namely the existence (assuming Szpiro's Conjecture) of an 
absolute constant $ C $ such that $ \tau < N^{C/\log\log N} $, but this was not needed in that paper. 
Actually I wanted $ \tau $ to be small, as the goal was to find curves  with exceptionally big 
Tate-Shafarevich groups, and the Birch and Swinnerton-Dyer Conjecture suggests that then $ \tau $ 
being small is helpful.

Recently this little auxiliary lemma has received interest from, a.o., Hector Pasten, in relation to 
the $ abc $ Conjecture, see \cite{P1}, \cite{P3}. In his paper \cite{P2} Pasten improved upon Lemma 1 
by the following conjecture.
\begin{conjecture}[Pasten, 2020]
For all elliptic curves 
\[ \tau \ll N^{\left(\frac73\log 3 + \epsilon\right)/\log\log N} . \]
\end{conjecture} 
Note that $ \frac73\log 3 = 2.563\ldots $, and note that I left out the term $ + O_{\epsilon}(1) $
from \cite[Conjecture 1.1]{P2}, as it seems superfluous. Pasten remarks that the constant 
$ \frac73\log 3 $ might not be optimal, but also believes there is `some evidence' for it. Anyway,
large $ \tau $'s are rare, as it is known \cite{GOT} that the average $ \tau $ is $ 1.8193\ldots $. 

This leads me to define the \emph{Tamagawa quality} of an elliptic curve by
\[ q_{\tau} = \dfrac{\log\tau\log\log N}{\log N} . \]
It is the purpose of this note to report on first experimental results searching for elliptic curves 
with an exceptionally high Tamagawa quality $ q_{\tau} $, and also for elliptic curves with an 
exceptionally large Tamagawa product $ \tau $ itself. The methods used here are restricted to picking 
the low hanging fruit only, and it is my hope that this note spurs interest from others to find 
better methods and results. Tables with all found curves can be found online \cite{dW3}.

\section{Curves from the LMFDB and Cremona Databases}
LMFDB \cite{LMFDB} is a database with web interface containing a wealth of data on a.o.\ elliptic 
curves. In particular it contains a collection, called 
\verb=ec_mwbsd=\footnote{\url{https://www.lmfdb.org/api/ec_mwbsd/}}, of data related to the Birch and
Swinnerton-Dyer Conjecture for $ 3824372 $ curves, containing for each curve a.o.\ the conductor and 
the Tamagawa product. So this is a good place to start. However, the web interface is not well suited
for directly checking the Tamagawa quality for all curves in the database. I found it easiest to
download the underlying database in two parts. 

The full data for all $ 3064705 $ curves with conductor up to $ 500000 $ can be downloaded directly
from John Cremona's database \cite{C}. I did so, and found the following results (see the green dots 
in Figure 1 below).
\begin{itemize}
\item $ 10795 $ curves have $ q_{\tau} > 1.5 $, 
\item $ 135 $ curves have $ q_{\tau} > 2 $, 
\item the curve with largest $ q_{\tau} = 2.30681 $ is 
      $ y^2 + x y = x^3 - 1054050116 x - 12046088636400 $, with  
      $ N = 39270 = 2 \, 3 \, 5 \, 7 \, 11 \, 17 $ and $ \tau = 31104 = 2^7 \, 3^5 $, 
\item the curve with largest $ \tau = 87040 = 2^{10} \, 5 \, 17 $ is
      $ y^2 + x y = x^3 - 4456595642213 x - 1538486355950810000 $, with  
      $ N = 364650 = 2 \, 3 \, 5^2 \, 11 \, 13 \, 17 $ and $ q_{\tau} = 2.26473 $.
\item Full tables \verb!output_cremona_qua.txt!, \verb!output_cremona_tam.txt! are on \cite{dW3}.
\end{itemize}
Note that Cremona's Database contains all curves with conductor up to $ 500000 $.

This leaves $ 759667 $ curves from the LMFDB \verb=ec_mwbsd= collection with conductor between 
$ 500000 $ and $ 300000000 $. The web interface does allow an easy download of data through
\url{https://www.lmfdb.org/EllipticCurve/Q/?conductor=500000-}, but this does not include Tamagawa 
products, one only gets conductors and Weierstrass coefficients $ a_1, a_2, a_3, a_4, a_6 $. However, 
computing $ \tau $ by SageMath \cite{S} then is trivial, with the code snippet \\
\verb!tau = EllipticCurve([a1,a2,a3,a4,a6]).tamagawa_product()! \\
Doing this I found the following results.
\begin{itemize}
\item no curve has $ q_{\tau} > 1.5 $, 
\item the largest $ q_{\tau} $ found is $ 1.22859 $, 
\item the largest $ \tau $ found is $ 576 $.
\end{itemize}
This somewhat disappointing result is probably due to the fact that of the curves of conductor
between $ 500000 $ and $ 300000000 $ only those of prime or $ 7 $-smooth conductor have been
incorporated, whereas $ \tau $ seems to get large only when the conductor has many small prime 
factors.

\section{Curves from $ abc $-triples}
\subsection{$ abc $-triples and Frey-Hellegouarch curves}
Let $ a, b, c $ be a triple of positive integers satisfying $ a + b = c $, $ a < b $, and 
$ \gcd(a,b) = 1 $. Its \emph{radical} $ r(a,b,c) $ is the product of the distinct prime divisors of 
$ a $, $ b $ and $ c $, i.e.\ $ r(a,b,c) = \displaystyle \prod_{\text{prime }p|abc} p $. Such a 
triple $ a, b, c $ is called an \emph{$ abc $-triple} if it satisfies $ c > r(a,b,c) $. The $ abc $ 
Conjecture states that $ c \ll r(a,b,c)^{1+\epsilon} $, and this has led to the definition of the
\emph{quality} of an $ abc $-triple as $ q(a,b,c) = \dfrac{\log c}{\log r(a,b,c)} $ (so a triple is 
an $ abc $-triple if and only if $ q(a,b,c) > 1 $). This concept of quality, although not yet under 
that name, seems to appear for the first time in my paper \cite{dW1}. Later also the \emph{merit} 
$ m(a,b,c) = (q(a,b,c)-1)^2 \log r(a,b,c) \log\log r(a,b,c) $ was introduced as an interesting
measure for $ abc $-triples. See Bart de Smit's website \cite{dS} for a wealth of experimental 
information on high quality and high merit $ abc $-triples. In this note I adopt the terminology 
\emph{medium quality} for an $ abc $-triple with $ 1.3 < q(a,b,c) < 1.4 $, next to \emph{high 
quality} for an $ abc $-triple with $ q(a,b,c) > 1.4 $ (as used in \cite{dW1}), and following 
\cite{dS} I also use \emph{high merit} for an $ abc $-triple with $ m(a,b,c) > 24 $, and 
\emph{unbeaten} if there is no $ abc $-triple known with larger $ c $ and larger $ q(a,b,c) $.

For an $ abc $-triple $ a, b, c $ it makes sense to look at its Frey-Hellegouarch curve, defined by
$ y^2 = x(x-a)(x+b) $, because its conductor equals $ r(a,b,c) $ up to a bounded power of $ 2 $, so
that such elliptic curves have exceptionally small conductors precisely for good quality
$ abc $-triples.

The Birch and Swinnerton-Dyer Conjecture for elliptic curves over $ \Q $ states
\[ \Omega \tau | \Sha | = \dfrac{T^2}{R} \lim_{s\to1} \dfrac{L(s)}{(s-1)^r} , \]
where $ \Omega = \omega $ or $ 2 \omega $ for the period $ \omega $, $ \Sha $ is the Tate-Shafarevich 
group, $ T $ is the order of the torsion group, $ R $ is the regulator, $ L(s) $ is the $ L $-series, 
and $ r $ is the rank of the elliptic curve. This conjecture is believed to hold for all elliptic
curves, but in this note I restrict to Frey-Hellegouarch curves.

The conductor is not explicitly there in the Birch and Swinnerton-Dyer formula, but its influence 
comes via $ \Omega $, as it is known that $ \Omega \ll \dfrac{\log c}{\sqrt{c}} $, which in the case 
of a good $ abc $-triple implies $ \Omega \ll N^{-1/2+\epsilon} $. In other words, the 
Frey-Hellegouarch curve then has an exceptionally small period, and the Birch and Swinnerton-Dyer 
Conjecture then suggests that this must be compensated for somewhere. The somewhat unusual way I have 
used above to present the Birch and Swinnerton-Dyer formula is suggestive for where I will be looking 
for this compensation.

The main idea of \cite{dW2} was that if one can show that, next to $ \Omega $, also the Tamagawa 
product $ \tau $ is small compared to the conductor, then one may expect large $ \Sha $'s, and this 
idea turned out to be fruitful, see also \cite{N}, \cite{DW}, \cite{B}. In this note I now 
complement this with the idea that, although there is a subpolynomial upper bound 
$ \tau \ll N^{C/\log\log N} $, it may still occur that large $ \tau $ accounts for a substantial part 
of this compensation of very small periods. In other words, one may expect big Tamagawa products also 
at Frey-Hellegouarch curves, and, like in \cite{dW2}, their quadratic twists and isogenous curves.

So this sets the program for the remainder of this note: to search for elliptic curves with large
Tamagawa products $ \tau $, and large Tamagawa qualities $ q_{\tau} $, by looking at curves isogenous
to quadratic twists of Frey-Hellegouarch curves for known good $ abc $-triples. For those 
$ abc $-triples the website of Bart de Smit \cite{dS} is an amazingly good source.

Let's start with a picture giving an overview of the results.

\begin{figure}[h] \centering
\captionsetup{justification=centering}
\includegraphics[width=\textwidth]{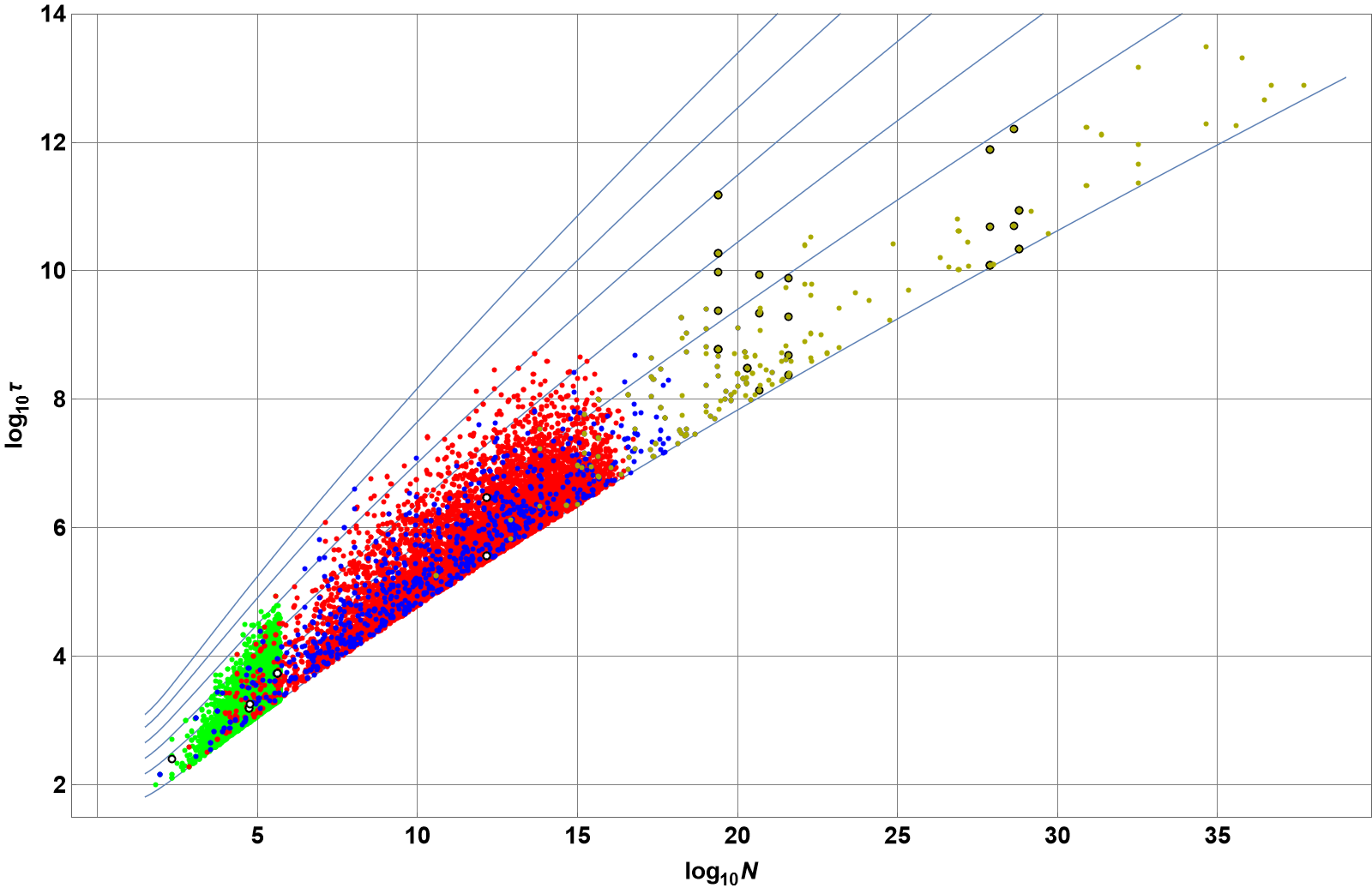}
\caption{Elliptic curves with large Tamagawa product $ \tau $ and Tamagawa quality 
$ q_{\tau} > 1.5 $. 
\\[1ex] {\small Legend: \begin{tabular}[t]{r@{\;}l}
curved lines: & $ q_{\tau} = \dfrac73\log 3 = 2.563\ldots, 2.4, 2.2, 2, 1.8, 1.5 $, \\
green dots: & curves from the Cremona and LMFDB databases, \\
yellow dots: & curves from high merit $ abc $-triples, \\
red dots: & curves from medium quality $ abc $-triples, \\
blue dots: & curves from high quality $ abc $-triples, \\
black circles: & curves from `triples from triples'.
\end{tabular}}}
\end{figure}

\subsection{Curves from high quality $ abc $-triples}
There are $ 241 $ high quality $ abc $-triples known, nicely presented as such on \cite{dS}. 
For each I took the twisted Frey-Hellegouarch curves $ d y^2 = x(x-a)(x+b) $ for 
$ d = \pm1, \pm2, \pm3, \pm5, \pm6 $, and some isogenous curves, computed by SageMath \cite{S} as 
follows: \\
\verb!E = EllipticCurve([0, d*(b-a), 0, -d^2*a*b, 0])! \\
\verb![E.isogeny(E(te)).codomain() for te in E.torsion_subgroup()]! \\
and then using SageMath's \verb!E.conductor()!, \verb!E.tamagawa_product()! to compute the conductor 
and the Tamagawa product (aborting the computation for each curve when it took more than 1 minute). 
This gave the following results, made visible in Figure 1 as blue dots.
\begin{itemize}
\item $ 841 $ curves have $ q_{\tau} > 1.5 $, 
\item $ 28 $ curves have $ q_{\tau} > 2 $, 
\item the curve with largest $ q_{\tau} = 2.39875 $ is \\
      $ y^2 + x y = x^3 - 2713479277841926834110 x - 53674762419393192464788215315900 $, \\
      with $ N = 105872910 = 2 \, 3 \, 5 \, 11 \, 13 \, 23 \, 29 \, 37 $ and 
      $ \tau = 3981312 = 2^{14} \, 3^5 $, \\
      isogenous to the Frey-Hellegouarch curve, twisted by $ -1 $, for the $ abc $-triple with \\
      $ a = 22771715409 = 3^{16} \, 23^2 $, \\
      $ b = 348972425216 = 2^{13} \, 29^2 \, 37^3 $, \\
      $ c = 371744140625 = 5^9 \, 11^4 \, 13 $, \\
      $ q(a,b,c) = 1.44181 $, $ m(a,b,c) = 10.5196 $, 
\item the curve with largest $ \tau = 152202903552 = 2^{28} \, 3^4 \, 7 $ is \\
      $ y^2 + x y = x^3 - 243293616838005191387643029131295594469691482466549330x - $ \\ 
      $ 46189598313302475345413359036293931548705009829803079259456070575837049845007100 $, \\
      with $ q_{\tau} = 2.18988 $ and $ N = 25180873035975641490 = 2 \, 3 \, 5 \, 11 \, 17 \, 19 \,
      23 \, 37 \, 43 \, 61 \, 127 \, 173 \, 4817 $, \\
      isogenous to the Frey-Hellegouarch curve, twisted by $ -1 $, for the $ abc $-triple with \\
      $ a = 44790692380548068359375 = 5^9 \, 17^2 \, 23^4 \, 37^2 \, 43 \, 4817 $, \\
      $ b = 3417300183328464869570036529 = 3^{14} \, 11^8 \, 61^2 \, 173^4 $, \\
      $ c = 3417344974020845417638395904 = 2^{52} \, 19^6 \, 127^2 $, \\
      $ q(a,b,c) = 1.41918 $, $ m(a,b,c) = 29.8237 $, \\
      and also for the $ abc $-triple with \\
      $ a = 146767394485224241 = 23^8 \, 37^4 $ \\
      $ b = 13669290314405085785446416384 = 2^{28} \, 3^7 \, 11^4 \, 19^3 \, 61 \, 127 \, 173^2 $ \\
      $ c = 13669290314551853179931640625 = 5^{18} \, 17^4 \, 43^2 \, 4817^2 $ \\
      $ q(a,b,c) = 1.45022 $, $ m(a,b,c) = 34.4028 $, twisted by $ -1 $; \\
      this peculiar situation is explained in Section 3.6.
\item Full tables \verb!output_high_quality_qua.txt!, \verb!output_high_quality_tam.txt! are on 
\cite{dW3}.
\end{itemize}

\subsection{Curves from medium quality $ abc $-triples}
There are $ 1947 $ medium quality $ abc $-triples known. They have been extracted from the two big
tables \verb!triples_below_1018_revised! (all $ 14482065 $ $ abc $-triples below $ 10^{18} $) and 
\verb!big_triples! (an additional $ 9345651 $ $ abc $-triples between $ 10^{18} $ and $ 2^{63} $).
With those medium quality $ abc $-triples I did exactly the same as I did with the high quality 
$ abc $-triples. This gave the following results, made visible in Figure 1 as red dots.
\begin{itemize}
\item $ 6342 $ curves have $ q_{\tau} > 1.5 $, 
\item $ 172 $ curves have $ q_{\tau} > 2 $, 
\item the curve with largest $ q_{\tau} = 2.39177 $ is \\
      $ y^2 + x y = x^3 - 986143769212695065 x - 376928045756312748465752775 $, \\
      with $ N = 13232310 = 2 \, 3 \, 5 \, 7 \, 13 \, 37 \, 131 $ and 
      $ \tau = 1228800 = 2^{14} \, 3 \, 5^2 $, \\
      isogenous to the Frey-Hellegouarch curve, twisted by $ -1 $, for the $ abc $-triple with \\
      $ a = 658489 = 13 \, 37^3 $, \\
      $ b = 6879707136 = 2^{20} \, 3^8 $, \\
      $ c = 6880365625 = 5^5 \, 7^5 \, 131 $, \\
      $ q(a,b,c) = 1.38137 $, $ m(a,b,c) = 6.67124 $,
\item the curves (isogenous) with largest $ \tau = 509607936 = 2^{21} \, 3^5 $ are \\
      $ y^2 + x y = x^3 - 13290632950903796089218578404113705 - $ \\ 
      $ 589748175639869043839018535079512741047463438442023 $, and \\
      $ y^2 + x y = x^3 - 13300785823649058521269530800913705 - $ \\ 
      $ 588802020491225519147911238670522467018762520682023 $, \\
      both with $ q_{\tau} = 2.19900 $ and $ N = 44947841915130 = 2 \, 3 \, 5 \, 7 \, 11 \, 17 \, 
      19 \, 23 \, 59 \, 103 \, 431 $, \\
      isogenous to the Frey-Hellegouarch curve, twisted by $ -1 $, for the $ abc $-triple with \\
      $ a = 40675641638471 = 7^3 \, 17^9 $, \\
      $ b = 798697622664921529 = 11^3 \, 19^3 \, 23 \, 59^2 \, 103^3 $, \\
      $ c = 798738298306560000 = 2^{20} \, 3^8 \, 5^4 \, 431^2 $, \\
      $ q(a,b,c) = 1.31127 $, $ m(a,b,c) = 10.5021 $. 
\item Full tables \verb!output_medium_quality_qua.txt!, \verb!output_medium_quality_tam.txt! are 
on \cite{dW3}.
\end{itemize}

\subsection{Curves from high merit $ abc $-triples}
There are $ 202 $ high merit $ abc $-triples available on \cite{dS}. For almost all I succeeded to
do the same procedure described above for high and medium quality $ abc $-triples. This produced the 
following results, made visible in Figure 1 as yellow dots.
\begin{itemize}
\item $ 190 $ curves have $ q_{\tau} > 1.5 $, 
\item $ 172 $ curves have $ q_{\tau} > 2 $, 
\item the curve with largest $ q_{\tau} = 2.18988 $ is the same as found with the largest $ \tau $ 
      for the curves coming from high quality $ abc $-triples, 
\item the curve with largest $ \tau = 30644423884800 = 2^{25} \, 3^4 \, 5^2 \, 11 \, 41 $ is \\
      $  y^2 + x y + y = x^3 - x^2 -
         621574482712904069167623787332097562003319547609729892578\backslash $ \\ 
      $ 282444666996972826570320570862 x - 5964690030130337213799773148403012416346828\backslash $ \\
      $ 0237789284970994254311306494753746869776095054625199343391902064488625755369\backslash $ \\
      $ 77569403651 $, \\
      with $ q_{\tau} = 1.70510 $ and $ N = 43081596887429422193675039055866970 = $ \\
      $ 2 \, 3^2 \, 5 \, 7 \, 11 \, 17 \, 19 \, 29 \, 43 \, 73 \, 83 \, 97 \, 103 \, 151 \, 577 \, 
      751 \, 3167 \, 1230379 $, \\
      isogenous to the Frey-Hellegouarch curve, twisted by $ -3 $, for the $ abc $-triple with \\
      $ a = 695606563606442148006101677581923 = 73^3 \, 97^2 \, 103^4 \, 577 \, 751 \, 3167 \,
      1230379 $, \\
      $ b = 57576591665034362126590541368210176398589952 = 2^{45} \, 17 \, 19^{10} \, 29^4 \, 43^5 \,
      151 $, \\
      $ c = 57576591665729968690196983516216278076171875 = 3^{11} \, 5^{33} \, 7^9 \, 11^2 \, 83^3 $, 
      \\
      $ q(a,b,c) = 1.28114 $, $ m(a,b,c) = 27.1356 $. 
\item Full tables \verb!output_high_merit_qua.txt!, \verb!output_high_merit_tam.txt! are on 
\cite{dW3}.
\end{itemize}

\subsection{Curves from unbeaten $ abc $-triples}
Finally, there are $ 160 $ unbeaten $ abc $-triples available on \cite{dS}, not necessarily different 
from curves in categories I have found above. Because those $ abc $-triples quickly get amazingly 
large, I was only able to process the 30 ones with smallest value of $ c $, and this did not yield
any examples not already found above in other categories. Full tables
\verb!output_unbeaten_qua.txt!, \verb!output_unbeaten_tam.txt! are on \cite{dW3}.

\subsection{Curves from triples from triples}
There is a nice trick to create new, hopefully better quality, $ abc $-triples from triples of lesser
quality, provided some miracle occurs. Assume $ a, b, c $ is an $ abc $-triple of reasonable quality,
and write $ d = a + c $, $ e = b + c $. Then look at the derived triples $ a, c, d $ and 
$ b, c, e $, and hope for the miracle that one of them is of reasonable quality as well, maybe even 
of quality $ > 1 $ so that it is an $ abc $-triple again. In that case I show how to get new triples
of probably better quality than $ q(a,b,c) $.

Note that $ a < b < c < d < e $, and observe that
\[ \left\{ \begin{array}{ll} c + a = d \\ c - a = b \end{array} \right. \Longleftrightarrow
   \left\{ \begin{array}{ll} d + b = 2c \\ d - b = 2a \end{array} \right. , \quad\quad
   \left\{ \begin{array}{ll} c + b = e \\ c - b = a \end{array} \right. \Longleftrightarrow
   \left\{ \begin{array}{ll} e + a = 2c \\ e - a = 2b \end{array} \right. . \]
Now put
\[ \begin{array}{lcll}
   (A_1,B_1,C_1) & = & (a^2,bd,c^2), & \\
   (A_2,B_2,C_2) & = & (b^2,4ac,d^2) & 
                       \text{ (if } b \text{ is even, divide by 4; if } b^2 > 4ac \text{, swap)}, \\
   (A_3,B_3,C_3) & = & (b^2,ae,c^2) & \text{ (if } b^2 > ae \text{, swap)}, \\
   (A_4,B_4,C_4) & = & (a^2,4bc,e^2) & \text{ (if } a \text{ is even, divide by 4)}. 
   \end{array} \]
Clearly all four new triples have $ A_i + B_i = C_i $ and $ \gcd(A_i,B_i) = 1 $ and $ A_i < B_i $, 
and to find out if they are $ abc $-triples only their quality has to be checked. Compared to the
$ abc $-triple $ a, b, c $, the numerator of the quality function almost doubles, namely from at most 
$ \log e < \log c + \log 2 $ to at least $ \log \frac14 c^2 = 2 \log c - 2 \log 2 $. The denominator 
however also increases, but most probably by a factor smaller than $ 2 $, as it grows from 
a three-term radical like $ r(a,b,c) $ to a four-term radical like $ r(a,b,c,d) $. In an ideal 
case where $ r(a) \approx r(b) \approx r(c) \approx r(d) \approx r(e) $ this growth factor in the 
denominator is $ \approx 4/3 $. So in such an ideal case the quality goes up by a factor 
$ \approx \dfrac{2}{4/3} = \dfrac{3}{2} $. No practical case is ideal, but I will give some examples
where the idea bears fruit, produces high-quality or high merit $ abc $-triples, which in turn
produce Frey-Hellegouarch curves with high Tamagawa quality.

Note that also $ a + e = 2c $ and $ b + d = 2c $, so that computing $ A_1, B_1, C_1 $ from 
$ a, 2b, e $  gives the same result as computing $ A_4, B_4, C_4 $ from $ a, b, c $, and computing 
$ A_2, B_2, C_2 $ from $ 2a, b, d $ gives the same result as computing $ A_3, B_3, C_3 $ from 
$ a, b, c $.

This idea of creating triples from triples is not new, it occurs in J.P.\ van der Horst's master 
thesis \cite{vdH}.

I tried this out for all $ abc $-triples found on \cite{dS}. This resulted in 13 examples with
quality above $ 1.4 $ or merit above $ 24 $. They are shown by black circles in Figure 1. Needless to 
say that they were all already present in these tables, so no new interesting curves with respect to 
Tamagawa quality were found, although the larger examples certainly lead to large $ \tau $ and
$ q_{\tau} $. But I found three cases of interest.

The first is from $ a = 10 = 2 \, 5 $, $ b = 2187 = 3^7 $, $ c = 2197 = 13^3 $, with 
$ q(a,b,c) = 1.28975 $, which for $ e = 4384 = 2^5 \, 137 $ leads to a not too bad quality 
$ q(b,c,e) = 0.90396 $, and then produces two high quality $ abc $-triples: \\
$ A_3 = 43840 = 2^6 \, 5 \, 137 $, $ B_3 = 4782969 = 3^{14} $, $ C_3 = 4826809 = 13^6 $, \\
with $ q(A_3,B_3,C_3) = 1.41370 $, and \\
$ A_4 = 25 = 5^2 $, $ B_4 = 4804839 = 3^7 \, 13^3 $, $ C_4 = 4804864 = 2^8 \, 137^2 $, \\
with $ q(A_4,B_4,C_4) = 1.41328 $. 

The second is a similar case, from $ a = 383102329 = 23^4 \, 37^2 $, 
$ b = 58457678566023 = 3^7 \, 11^4 \, 61 \, 173^2 $, $ c = 58458061668352 = 2^{26} \, 19^3 \, 127 $,
with $ q(a,b,c) = 1.13257 $, which for $ e = 116915740234375 = 5^9 \, 7^2 \, 43 \, 4817 $ leads to a 
not too bad quality $ q(b,c,e) = 0.854092 $, and then produces two high quality and high merit 
$ abc $-triples: \\
$ A_3 = 44790692380548068359375 = 5^9 \, 17^2 \, 23^4 \, 37^2 \,43 \, 4817 $, \\
$ B_3 = 3417300183328464869570036529 = 3^{14} \, 11^8 \, 61^2 \, 173^4 $, \\
$ C_3 = 3417344974020845417638395904 = 2^{52} \, 19^6 \, 127^2 $, \\
with $ q(A_3,B_3,C_3) = 1.41918 $, $ m(A_3,B_3,C_3) = 29.8237 $, and \\
$ A_4 = 146767394485224241 = 23^8 \, 37^4 $, \\
$ B_4 = 13669290314405085785446416384 = 2^{28} \, 3^7 \, 11^4 \, 19^3 \, 61 \, 127 \, 173^2 $, \\
$ C_4 = 13669290314551853179931640625 = 5^{18} \, 17^4 \, 43^2 \, 4817^2 $, \\
with $ q(A_4,B_4,C_4) = 1.45022 $, $ m(A_4,B_4,C_4) = 34.4028 $. \\
Those two $ abc $-triples appeared above, at the curve with the largest Tamagawa product found from 
high quality $ abc $-triples as well as the largest Tamagawa quality found from high merit 
$ abc $-triples.

The third is again a similar case, from $ a = 158810997195450625 = 5^4 \, 53^2 \, 67^6 $, \\
$ b = 4025783917396764928 = 2^8\, 23^2 \, 61^6 \, 577 $, 
$ c = 4184594914592215553 = 17^6 \, 311 \, 823^3 $,
with $ q(a,b,c) = 1.08916 $, which for 
$ d = 4343405911787666178 = 2 \, 3^4 \, 7 \, 13 \, 109^4 \, 2087219 $ leads to a 
not too bad quality $ q(a,c,d) = 0.847863 $, and then produces two high merit $ abc $-triples: \\
$ A_1 = 25220932830213426279247196812890625 = 5^8 \, 53^4 \, 67^{12} $, \\
$ B_1 = 17485613666400818352222466948402205184 = 2^9 \, 3^4 \, 7 \, 13 \, 23^2 \, 61^6 \, 109^4 \,
577 \, 20872194 $, \\
$ C_1 = 17510834599231031778501714145215095809 = 17^{12} \, 311^2 \, 823^6 $, \\
with $ m(A_1,B_1,C_1) = 30.0604 $, and \\
$ A_2 = 664559691245401291839706490968570625 = 5^4 \, 17^6 \, 53^2 \, 67^6 \, 311 \, 823^3 $, \\
$ B_2 = 4051734037392610655275387090022711296 = 2^{14} \, 23^4 \, 61^{12} \, 577^2 $, \\
$ C_2 = 4716293728638011947115093580991281921 = 3^8 \, 7^2 \, 13^2 \, 109^8 \, 20872194^2 $, \\
with $ m(A_2,B_2,C_2) = 26.5098 $. 

Full tables \verb!output_triples_from_triples_qua.txt!, 
\verb!output_triples_from_triples_tam.txt! are on \cite{dW3}.

\section{Discussion}
\begin{table}[h]
$ \begin{array}{|l|rr|r|r|r|} \hline
  \text{database} & \text{\# curves} & \text{\# $abc$-triples} & \text{\# curves} &
  \text{largest } q_{\tau} & \text{largest } \tau \\
  & & & \text{with } q_{\tau} > 1.5 & & \\ \hline
  \text{Cremona}           & 3064705 & & 10795 & 2.30681 &          87040 \\
  \text{LMFDB}             &  759667 & &     0 & 1.22859 &            576 \\
  \text{high quality}         & &  241 &   841 & 2.39875 &   152202903552 \\
  \text{medium quality}       & & 1947 &  6342 & 2.39177 &      509607936 \\
  \text{high merit}           & &  202 &   190 & 2.18988 & 30644423884800 \\
  \text{unbeaten}             & &  160 &    24 & 2.18988 & 20615843020800 \\
  \text{triples from triples} & &   13 &    42 & 2.18988 &  1594506608640 \\ \hline
  \text{all together}         & &      & 17760 & 2.39875 & 30644423884800 \\ \hline
  \end{array} $
\caption{Overview of found data on curves from different sources. Full data on \cite{dW3}.}
\end{table}

For a summary of my results, see Table 1, and also Figure 1. Note that not all collections are 
necessarily disjunct.

As a first conclusion it seems justified to state that Frey-Hellegouarch curves for good 
$ abc $-triples are a good source for high quality Tamagawa products. Finding large $ \tau $ is 
probably not that hard, just by looking at curves with a large number of bad primes; the curve with 
largest $ \tau $ found here ($ \tau \approx 3 \times 10^{13} $) is a good example of that, having 
$ 18 $ bad primes. But finding curves with a high Tamagawa quality is another matter. 

One might describe Pasten's conjecture 2 as cited above by the simple inequality
\[ \limsup_{E/\Q} q_{\tau} \leq \dfrac{7}{3}\log 3 . \]
Pasten already remarks that it could be open for improvement, and my experiments so far seem not to 
enthousiastically support a conjecture that $ \dfrac{7}{3}\log 3 $ is the true constant. As may be 
clear I am also interested in a lower bound for $ \displaystyle\limsup_{E/\Q} q_{\tau} $. My 
experiments do not shed much light on this problem. It seems plausible to me to conjecture that at 
the very least $ \displaystyle\limsup_{E/\Q} q_{\tau} > 0 $, and I leave it to others to elaborate, 
and to come up with an argumented conjectured value, or at least a positive lower bound, for 
$ \displaystyle\limsup_{E/\Q} q_{\tau} $. It does seem that Frey-Hellegouarch curves show split 
multiplicative reduction at a substantial portion of the bad primes, causing the local Tamagawa 
numbers at such primes being the exponents of the primes in the discriminant $ 16 (abc)^2 $ and thus 
contributing to large $ \tau $'s, and this gives some hope that Pasten's analysis points in the right 
direction.

I also leave it for others to investigate whether the methods of papers searching for large 
$ \Sha $'s, in particular \cite{N}, \cite{DW} and \cite{B}, now geared towards small $ \tau $ for
obvious reasons, can be adapted to searching for big $ \tau $'s at not too big conductors as well.

The code I wrote for this little project is so embarrassingly trivial that I did not want to set up 
an online repository for it. Some of the most essential code snippets are mentioned in the text
above; everything else is simple data manipulation.

\end{document}